\documentclass[11pt,bezier,epsf]{amsart}
\usepackage{amsmath,amssymb,amsfonts}
\usepackage{graphicx}
\usepackage{cite}
\usepackage{subfigure}
\newbox\subfigbox
\makeatletter
{\def\caption##1{\gdef\subcapsave{\relax##1}}%
  \let\subcapsave\@empty%
  \setbox\subfigbox\hbox%
  \bgroup}%
{\egroup%
  \subfigure[\subcapsave]{\box\subfigbox}}%
\makeatother

 \newtheorem{thm}{Theorem}[section]
 \newtheorem{cor}[thm]{Corollary}

 \newtheorem{prop}[thm]{Proposition}
 
\theoremstyle{definition}
\newtheorem{example}[thm]{Example}
 \newtheorem{defn}[thm]{Definition}
\theoremstyle{remark}
 \newtheorem{rem}[thm]{Remark}
\numberwithin{equation}{section}

\begin{document}

\title[C*-algebraic approach to fixed point theory]{C*-algebraic approach to fixed point theory generalizes Baggett's theorem to groups with discrete reduced duals}

\author[F. Naderi]{Fouad Naderi}

\address{Department of Mathematical and Statistical Sciences, University of
Alberta, Edmonton, Alberta, T6G 2G1, Canada.}
\email{naderi@ualberta.ca}


\begin{abstract}
In this paper, we show that if the reduced Fourier-Stieltjes algebra $B_{\rho}(G)$ of a second countable locally compact group $G$ has either weak* fixed point property
or asymptotic center property, then $G$ is compact. As a result, we give affirmative
answers to open problems raised by Fendler and et al. in 2013. We then conclude that a second countable group with a discrete reduced dual must be compact. This generalizes a theorem of Baggett. We also construct a compact scattered Hausdorff space $\Omega$ for which the dual of the scattered C*-algebra $C(\Omega)$ lacks weak* fixed point property. The C*-algebra $C(\Omega)$ provides a negative answer to a question of Randrianantoanina in 2010. In addition, we prove a variant of Bruck's generalized fixed point theorem for the preduals of von Neumann algebras. Furthermore, we give some examples emphasizing that the conditions in Bruck's generalized conjecture (BGC) can not be weakened any more.
\end{abstract}

\maketitle
\section{Introduction}
Let $K$ be a subset of a Banach space $E$. A self mapping $T$ on $K$ is said to
be {\it non-expansive} if  $\|T(x)-T(y)\|\leq \|x-y\|$ for all $x,y\in K$. We
say that $E$ has the {\it fixed point property (fpp)} if for every bounded
closed convex non-empty subset $K$ of $E$, any non-expansive self mapping on $K$
has a fixed point. $E$ is said to have the {\it weak fixed point property (weak
fpp)} if for every weakly compact convex non-empty subset $K$ of $E$, any
non-expansive self mapping on $K$ has a fixed point. For a dual Banach space
$E$, {\it weak* fixed point property (weak* fpp)} is defined similarly. One can
easily check that the weak fpp is the weakest notion among these properties, this is
why the weak fpp is the most desirable property among them. Note that, the definition of
fpp is based on bounded closed sets not on compact sets, otherwise; the fpp
would be automatic by Schauder's fixed point theorem.

In recent years, Lau has launched a program in which he has successfully related
various topological properties of a locally compact group $G$ to fixed point
properties of Banach algebras associated to $G$ ( see \cite{Fendler-Lau},
\cite{Lau-Leinert}, \cite{Lau-Mah-Ulger},  \cite{Lau-Ulger}, and
\cite{Randrian}). In \cite{Randrian}, it has been shown that the
Fourier-Stieltjes algebra $B(G)$ has weak fpp if and only if $G$ is an
[AU]-group. Fendler and et al. in \cite{Fendler-Lau} have shown that $B(G)$ has
weak* fpp if and only if $G$ is a compact group. They also have shown that these
are equivalent to saying that $B(G)$ has asymptotic center property (see
Definition \ref{asym}). Then, they asked if the corresponding result is true for the reduced
Fourier-Stieltjes algebra $B_{\rho} (G)$ or not? In this paper we are mainly concerned with this
question. It turns out that the affirmative answer to this question has many important implications.

This paper is organized as follows: section 2 gives some preliminaries. In section 3,
we prove that the conjecture of Fendler
and et al. is true for second countable locally compact groups. We then show a second countable group with a discrete reduced dual is compact. This is a generalization of Baggett's theorem \cite{Bagget}. Also, our proof gives an alternative proof to Baggett's original theorem for a second countable locally compact group with a discrete full dual. Besides, for the compact scattered set $\Omega=\{\frac{1}{n}: n\in \mathbb{N} \}\cup \{0\} $ in real numbers, we show that the dual space of the C*-algebra $C(\Omega)$ does not have weak* fpp. This gives a negative answer to \cite[Question 3.9]{Randrian}.
In section 4, we discuss
about Bruck's generalized fixed point theorem. Here, we show that if the predual of a von Neumann algebra has weak fpp, then it must have weak fpp for left reversible semigroups. Also, we give examples of
semigroups acting separately continuous on compact sets such that every $s\in S$ has a
fixed point, but $S$ still does not have a common fixed point. The problem with these examples is that
they don't enjoy one or two conditions of Bruck's generalized conjecture. The paper concludes
with some open problem arising from the current work.

\section {Some preliminaries}
Let $G$ be locally compact group with a fixed left Haar measure $\mu$. Let
$L^{1}(G)$ be the group algebra of $G$ with convolution product. We define
$C^{\ast}(G)$, the group $C^{\ast}$-algebra of $G$, to be the completion of
$L^{1}(G)$ with respect to the norm
 $$\|f\|_{*}=sup\|\pi_{f}\|$$
where the supremum is taken over all non-degenerate $\ast$-representation $\pi$
of $L^{1}(G)$ as a $\ast$-algebra of bounded operators on a Hilbert space. Let
$\mathcal{B}(L^{2}(G))$ be the set of all bounded operators on the Hilbert space
$L^{2}(G)$ and $\rho$ be the left regular representation of $G$, i.e., for each
$f\in L^{1}(G)$, $\rho(f)$ is the bounded operator in $\mathcal{B}(L^{2}(G))$
defined by $\rho(f)(h)=f\ast h$, the convolution of $f$ and $h$ in $L^{2}(G)$.
The reduced group C*-algebra ${{C^{*}}}_{\rho}(G)$ is the closure of $\rho(L^{1}(G))$ 
in the norm topology of $\mathcal{B}(L^{2}(G))$ while the von Neumann algebra generated by $\rho$
, $VN(G)$, is
the closure of the same set in the weak operator topology. For any other unitary representation $\pi$, the von Neumann algebra generated by $\pi$
is the closure of $\pi(L^{1}(G))$ 
in the weak operator topology of $\mathcal{B}({\mathcal{H}}_{\pi})$, where ${\mathcal{H}}_{\pi}$ is the Hilbert space 
of the representation $\pi$.
Denote the set of all continuous positive definite function on $G$ by $P(G)$, and
the set of all continuous function on $G$ with compact support by $C_{c}(G)$. We
define the Fourier-Stieltjes  algebra of $G$, denoted by $B(G)$, to be the
linear span of $P(G)$. Then, $B(G)$ is a Banach algebra with pointwise
multiplication and the norm of each $\phi \in B(G)$ is defined by
$$\|\phi\|=sup\{ | \int f(t)\phi(t) d\mu(t) | : f\in L^{1}(G), \|f\|_{*}\leq 1
\}.$$
The Fourier algebra of $G$, denoted by $A(G)$, is defined to be the closed
linear span of $P(G)\cap C_{c}(G)$. Note that, $A(G)$ is the set of all matrix
elements of $\rho$ while $B(G)$ is the set of all matrix elements of all
continuous unitary representations of $G$. Clearly, $A(G)=B(G)$ when $G$ is
compact. It is known that $B(G)=C^{\ast}(G)^{\ast}$, where the duality is given
by $\langle f,\phi\rangle = \int f(t)\phi(t) d\mu(t), f\in L^{1}(G), \phi \in
B(G)$. Besides, we can define the reduced Fourier-Stieltjes algebra $ B_{\rho}
(G)$ of $G$ to be the dual space of the reduced group C*-algebra ${C^{*}}_{\rho}
(G)$. The duality is pretty much the same as $B(G)=C^{\ast}(G)^{\ast}$. When $G$
is amenable (e.g., abelian and compact groups), one has ${C^{*}}_{\rho}
(G)=C^{*}(G)$ and $B_{\rho} (G)=B(G)$. Furthermore, $A(G)$ and $B_{\rho}(G)$ are
closed ideals in $B(G)$ and $A(G)\subseteq B_{\rho} (G)\subseteq B(G)$ (see \cite{Eymard} for details).

The set of all unitarily equivalent classes of
all irreducible representations of a C*-algebra $\mathcal{A}$ is called the
spectrum of $\mathcal{A}$ and is denoted by $\hat{\mathcal{A}}$. The dual
of a locally compact group $G$, denoted by $\hat{G}$, is the spectrum of 
the group C*-algebra $C^{*}(G)$. Similarly, the reduced
dual of a locally compact group $G$, denoted by $\hat{G}_{\rho}$, is the
spectrum of the reduced C*-algebra ${{C^{*}}}_{\rho}(G)$. The reduced dual has
many equivalent definitions (see \cite{Dixmier}).

A von Neumann algebra $\mathcal{M}$ is called {\it atomic} if every nonzero
(self-adjoint) projection in $\mathcal{M}$ dominates a minimal nonzero projection
in $\mathcal{M}$. Let $G$ be a locally compact group. $G$ is called an
[AU]-group if the von Neumann algebra generated by any continuous unitary
representation of $G$ is atomic, i.e., every continuous unitary representation
of $G$ is completely decomposable. $G$ is called an [AR]-group if the von Neumann
algebra generated by $\rho$, $VN(G)$, is atomic. Evidently, every [AU]-group is
an [AR]-group. $G$ is called an [IN]-group if there is a compact neighborhood of
the identity in $G$ which is invariant under inner automorphisms. For
these classes of groups, we have $[AU]\cap [IN]=[compact]$ (see \cite{Taylor} for more
details). $G$ is called a [Moore]-group, if each irreducible unitary
representation of $G$ is finite dimensional. For example, any compact group is a
Moore group \cite[Chapter
5]{Folland}. The converse is false since abelian groups are also
[Moore]-groups and can be non-compact \cite[Chapter 4]{Folland}.

 A Banach space $E$ is said to have {\it Dunford-Pettis property (DPP),}
 if for any Banach space $F$, every weakly compact linear operator
$u:E\longrightarrow F$ sends weakly Cauchy sequences into norm convergent
sequences. $E$ is said to have the {\it Radon-Nikodym property
(RNP)} if each closed convex sub-set $D$ of $E$ is dentable, i.e., for any
$\varepsilon >0$ there exists an $x$ in $D$ such that $ x\notin \overline{co}(D
\setminus B_{\varepsilon}(x))$, where $B_{\varepsilon}(x)=\{y\in X:
\|x-y\|<\varepsilon \}$ and $\overline{co}(K)$ is the closed convex hull of a
set $K\subseteq E$. The book \cite{Dis} contains a complete account of these properties.
A dual Banach space $E$ is said to have the {\it  uniform weak* Kadec-Klee property (UKK*)} if
for every $\epsilon >0$ there is $0<\delta<1$ such that for any subset $C$ of its closed unit ball
and any sequence $(x_n)$ in $C$ with $sep(x_n):=inf\{\rVert\ {x_n} -{ x_m}\rVert: n\neq m\}>\epsilon$, there is an $x$ in the weak* closure of $C$ with $\rVert x \rVert<\delta$. This property was firs introduced by van Dulst and Sims \cite{Dulst}. They proved that if
E has property UKK*, then $E$ has the weak* fpp. A dual Banach space $E$ is said to have the {\it weak* Kadec-Klee property (KK*)} if the weak* and norm convergence for sequences coincide on the unit sphere of $E$. It is well known that UKK* implies KK* ( see \cite{Lennard} and \cite{Lau-Mah 88} for more details).

Let $S$ be a {\it semi-topological semigroup}, i.e., $S$ is a semigroup with a
Hausdorff topology such that for each $a \in S$, the mappings $s\mapsto sa$ and
$s\mapsto as$ from $S$ into $S$ are continuous. $S$ is called {\it left (right)
reversible} if any two closed right (left) ideals of $S$ have non-void
intersection.

An {\it action} of $S$ on a subset $K$ of a topological space $E$ is a mapping
$(s,x)\mapsto s(x)$ from $ S \times K$ into $K$ such that $(st)(x)=s(t(x))$ for
$ s,t\in S, x\in K$. The action is {\it separately continuous} if it is
continuous in each variable when the other is kept fixed. When we speak of
separately continuous actions, we always put the norm topology on $K$ although
$K$ can be a weak or weak* compact set. Every action of $S$ on $K$ induces a
representation  of $S$ as a semigroup of self-mappings on $K$ denoted by
$\mathcal{S}$, and the two semigroups are usually identified. When the action is
separately continuous, each member of $\mathcal{S}$ is continuous. We say that
$S$ has a {\it common fixed point} in $K$ if there exists a point $x$ in $K$
such that $sx=x$ for all $s\in S$. When $E$ is a normed space, the action of $S$
on $K$ is {\it non-expansive} if $ \| s(x)-s(y)\| \leq \| x-y\|$ for all  $s \in
S$ and $x,y \in K$. There are also other types of action for a semi-topological
semigroup (see \cite{Holmes} and \cite{Amini}).

We say that a Banach space $E$ has the {\it fpp for left
reversible semigroups} if for every bounded closed convex non-empty subset $K$
of $E$, any non-expansively separately continuous action of a left reversible
semi-topological semigroup $S$  on $K$ has a fixed point. $E$ is said to have the
{\it weak fpp for left reversible semigroups} if for every
weakly compact convex non-empty subset $K$ of $E$, any non-expansively
separately continuous action of a semi-topological semigroup $S$  on $K$ has a
fixed point. For a dual Banach space $E$, {\it weak* fpp for
left reversible semigroups} is defined similarly.

Let $\mathbb{T}$ be the unit circle in the complex plane which is a locally
compact group under usual multiplication and Euclidean norm. Suppose that
$\mathbb{Z}$ is the group of all integers under addition and discrete topology.
By duality theory for abelian groups, one has $\hat{\mathbb{Z}}=\mathbb{T}$ and
$\hat{\mathbb{T}}=\mathbb{Z}$ (see \cite[Chapter 4]{Folland}). Alspach
\cite{Alspach} has shown that $L^{1}[0,1]=L^{1}(\mathbb{T})=A(\mathbb{Z})$ does
not have weak fpp (for abelian semigroups). This
inspired Lau and Mah in \cite{Lau-Mah} to show that for an [IN]-group $G$,
$A(G)$ has weak fpp for left reversible semigroups if and only
if $G$ is compact. Later, Randrianantoanina \cite{Randrian} showed that $A(G)$ (
resp. $B(G)$) has weak fpp for left reversible semigroups if
and only if $G$ is an [AR]-group ( resp. [AU]-group). Also, Lim's fixed point
theorem \cite{Lim} for dual Banach spaces shows that
$l^{1}(\mathbb{Z})=B(\mathbb{T})$ has weak* fpp for left
reversible semigroups. This inspired Lau and et al. in \cite{Fendler-Lau} and
\cite{Lau-Mah} to show that $B(G)$ has weak* fpp for left
reversible semigroups if and only if $G$ is compact.


\section {Weak* fixed point property of reduced Fourier-Stieltjes algebra and Baggett's generalized theorem}

The following theorem answers the first problem in
\cite[p.300]{Fendler-Lau} positively. Before stating it, we would like to thank Professor Gero Fendler for finding out a gap in our original proof and suggesting the idea of decomposing a separable C*-algebra with discrete spectrum according to \cite[10.10.6]{Dixmier}.

\begin{thm} \label{main.thm}
If the reduced Fourier-Stieltjes algebra $B_{\rho} (G)$ of a second countable locally
compact group has the weak* fpp, then $G$ is compact.
\end{thm}

{\bf Proof.} For a second countable locally compact group $G$, the reduced group C*-algebra  ${{C^{*}}}_{\rho}(G)$ is separable. Hence, the unit ball of $B_{\rho}(G)$ as the dual Banach space of the separable C*-algebra ${{C^{*}}}_{\rho}(G)$ is metrizable in weak* topology \cite[p.134]{Conway}. According to \cite[Theorem 2.1]{Fendler-Leinert}, the reduced dual $\hat{G}_{\rho}$ must be discrete. By the separability assumption and \cite[10.10.6]{Dixmier}, one gets the following decompositions
\begin{center}
${{C^{*}}}_{\rho}(G)=c_{0}-\oplus_{i\in I} \mathcal{K}({\mathcal{H}}_{i})$; \quad
$B_{\rho}(G)=l^{1}-\oplus_{i\in I} \mathcal{T}({\mathcal{H}}_{i})$
\end{center}
where each $\mathcal{H}_{i}$ is a separable Hilbert space, and  $\mathcal{K}({\mathcal{H}}_{i})$ and $\mathcal{T}({\mathcal{H}}_{i})$ are compact and trace class operators on it. By a digonalization argument,
one can embed  $B_{\rho}(G)$ into 
\begin{center}
$\mathfrak{T}=\mathcal{T}(l^{2}-\oplus_{i\in I}{\mathcal{H}}_{i})\cong {\mathcal{K}(l^{2}-\oplus_{i\in I}{\mathcal{H}}_{i})}^{*}={\mathfrak{K}}^{*}$
\end{center}

By \cite{Lennard}, $\mathfrak{T}$ has UKK*. Hence, it has KK*. The space $\mathfrak{K}$ is separable, so the unit ball of $\mathfrak{T}$ is metrizable in the weak* topology stemming from $\mathfrak{K}$. This alongside with the property KK* means that the norm and weak* topology on unit sphere of $\mathfrak{T}$ agree. On the other hand, the weak* topology of $\mathfrak{T}$ relativized to $B_{\rho}(G)$ comes from ${{C^{*}}}_{\rho}(G)$. So, we are allowed to say that the norm and weak* topology on unit sphere of $B_{\rho}(G)$ agree, as well. By applying \cite[Theorem 4.2]{Miao}, we deduce that $G$ is compact. $\blacksquare$

The following theorem is a consequence of Theorem \ref{main.thm} and generalizes a result of Baggett \cite{Bagget}. The proof may seem simple, but it is only after a great deal of theory that we can give such a proof.
\begin{thm} \label{Bagget}
A second countable group having a discrete reduced dual is compact.
\end{thm}
{\bf Proof.} Let $G$ be a second countable group and $\hat{G}_{\rho}$ be its discrete reduced dual. The spectrum of the reduced C*-algebra ${C^{*}}_{\rho}(G)$ is $\hat{G}_{\rho}$. According to \cite{Fendler-Leinert}, $ B_{\rho} (G)$ has the weak* fpp. By Theorem \ref{main.thm}, $G$ is compact. $\blacksquare$

\begin{rem}
The method of the proofs of Theorems \ref{main.thm} and \ref{Bagget} can be used to give an alternative proof to Baggett's theorem \cite{Bagget} regarding the full dual.
\end{rem}

We can also prove the converse of Theorem \ref{main.thm}. Although it can be
inferred from \cite[Theorem 4.7]{Fendler-Lau}, we prefer a direct proof to show the sophisticated ideas of Randrianantoanina \cite{Randrian} in studying fixed point theory. We use a well-known fact that  ${{C^{*}}}_{\rho}(G)=C^{*}(G)$ for every amenable group $G$.  The discussion given here can be viewed as a simplification of the proof of \cite[Theorem 4.2]{Lau-Mah} in light of von Neumann algebraic techniques of  \cite{Randrian}. For the sake of completeness we give some details from \cite{Lau-Mah}.

\begin{prop} \label{converse}
	If $G$ is compact, then the reduced Fourier-Stieltjes algebra $ B_{\rho} (G)$
	has weak* fpp for left reversible semi groups.
\end{prop}
{\bf Proof.} When $G$ is compact, $G$ is amenable and every element of $\hat{G}$ is finite
dimensional \cite[Chapter 5]{Folland}. So, $G$ is a Moore group and by
\cite[Theorem 4.2]{Lau-Ulger}, $C^{*}(G)={{C^{*}}}_{\rho}(G)$ has DPP. According
to \cite[Lemma 4.1]{Lau-Ulger}, $ {{C^{*}}}_{\rho}(G)$ must be a $c_0$-direct
sum of finite dimensional C*-algebras. But each finite dimensional C*-algebra,
is a direct sum of matrix algebras \cite[p.194]{Murphy}, so
\begin{center}
	${{C^{*}}}_{\rho}(G)=c_{0}-\oplus_{i\in I}\mathcal{K}({\mathcal{H}}_{i})$ and
	$B_{\rho}(G)=l^{1}-\oplus_{i\in I}\mathcal{T}({\mathcal{H}}_{i})$
\end{center}
where each ${\mathcal{H}}_{i}$ is a finite dimensional Hilbert space and
$\mathcal{K}({\mathcal{H}}_{i})$ and $\mathcal{T}({\mathcal{H}}_{i})$ denote
compact and trace class operators on ${\mathcal{H}}_{i}$. By \cite[Corrolary
3.7]{Randrian}, it is now obvious that $B_{\rho}(G)=l^{1}-\oplus_{i\in
I}\mathcal{T}({\mathcal{H}}_{i})$ has  the weak* fpp for left
reversible semi groups.$\blacksquare$

\begin{defn} \label{asym}
(a) Let $C$ be a non-empty subset of a Banach space $E$ and $\{ W_{\alpha} :
\alpha \in A\}$ be a decreasing net of non-empty bounded subsets of $E$. For
each $c\in C, \alpha \in A$ define
$$ r_{\alpha}(c)=sup\{ \parallel w-c\parallel: w\in W_{\alpha} \} $$
$$r(c)=inf\{r_{\alpha}(c) : \alpha \in A \}$$
$$r=inf\{ r(c) : c\in C \}$$

The set $ AC(\{W_{\alpha} :\alpha \in A\},C)=\{ c\in C : r(c)=r \}$ (the number
$r$ ) will be called the {\it asymptotic center (asymptotic radius)} of $\{
W_{\alpha} : \alpha \in A\}$ in $C$. This notion is due to Edelstein
\cite{Edelstein}. A dual Banach space $E$ is said to have the {\it asymptotic center property}
if for any non-empty weak* closed convex subset $C$ in $E$ and any decreasing
net $\{ W_{\alpha} : \alpha \in A\}$ of non-empty bounded subsets of $C$, the
set $ AC(\{W_{\alpha} :\alpha \in A\},C)=\{ c\in C : r(c)=r \}$ is a non-empty
norm compact convex subset of $C$. This property was first used by Lim
\cite{Lim}.

(b) A dual Banach space $E$ is said to have the {\it lim-sup property} if for any
decreasing net $\{ W_{\alpha} : \alpha \in A\}$ of non-empty bounded subsets of
$E$, and any weak* convergent net ${(\varphi_{\mu})}_{\mu \in M}$ with weak*
limit $\varphi$, we have:
$$ lim sup_{\mu} \|{\varphi}_{\mu}-\varphi\|+ r(\varphi)=lim sup_{\mu}
r({\varphi}_{\mu}) $$
, i.e.,
$$lim sup_{\mu} \|{\varphi}_{\mu}-\varphi\|+ lim sup_{\alpha} \{
\|\varphi-\psi\| : \psi \in W_{\alpha} \} =lim sup_{\mu}lim sup_{\alpha} \{
\|{\varphi}_{\mu}-\psi\| : \psi \in W_{\alpha} \}.$$

This notion was first used by Lim \cite{Lim}. Later, Lau and et al.
generalized it to arbitrary dual Banach spaces (see \cite[Lemma 3.1]{Lau-Mah}
and \cite[Lemma 5.1 and Defenition 5.2]{Fendler-Lau}).
\end{defn}

The next corollary answers the second question in \cite[p.300]{Fendler-Lau}.

\begin{cor} \label{asym.cpt}
If the reduced Fourier-Stieltjes algebra $ B_{\rho} (G)$ of a second countable locally
compact group has asymptotic center property, then $G$ is compact.
\end{cor}

{\bf Proof.} By \cite[Theorem 6.1, conditions b and c]{Fendler-Lau}, $ B_{\rho}
(G)$ has the weak* fpp. Now, apply Theorem
\ref{main.thm}.$\blacksquare$

The following theorem characterize the compactness of $G$ in terms of weak* fpp and
geometric properties of $ B_{\rho}(G)$. It is the reduced counterpart of \cite[Theorem
5.3]{Fendler-Lau}. Recall that a C*-algebra $\mathcal{A}$ is called {\it scattered} if every positive linear
functional on $\mathcal{A}$ is the sum of a sequence of pure states.

\begin{thm}
Let $G$ be a second countable locally compact group. The followings are equivalent:

(a) $G$ is compact.

(b) $B_{\rho} (G)$  has the lim-sup property.

(c) $B_{\rho} (G)$  has the asymptotic center property.

(d) $ B_{\rho} (G)$ has the weak* fpp for left reversible
semi-topological semigroups.

(e) $B_{\rho}(G)$ has the weak* fpp for non-expansive
mappings.

(f) $lim sup_{\mu} \|{\varphi}_{\mu}-\varphi\| + \|\varphi\|=lim sup_{\mu}
\|{\varphi}_{\mu}\|$ for any bounded net $(\varphi_{\mu})$ in $ B_{\rho} (G)$
which converges to $\varphi \in B_{\rho} (G)$ in the weak* topology.

(g) for any net $(\varphi_{\mu})$ and any $\varphi \in B_{\rho} (G)$ we have
that $\|\varphi_{\mu}-\varphi\|\rightarrow 0$ if and only if
$\varphi_{\mu}\rightarrow \varphi$ in the weak* topology and
$\|\varphi_{\mu}\|\rightarrow \|\varphi\|$

(h) On the unit sphere of $ B_{\rho} (G)$ the weak* and norm topology are the
same.

(i) Each member of $\hat{G}_{\rho}$ is finite dimensional and $B_{\rho}(G)$ has RNP.

(j) $B_{\rho}(G)$ has RNP and DPP.

(k) The C*-algebra ${C^{*}}_{\rho} (G)$ is scattered and has DPP.

(l) The reduced dual $\hat{G}_{\rho}$ is discrete.

\end{thm}
{\bf Proof.}

(a)$\Rightarrow$(b): When $G$ is compact, $G$ is amenable, and then
$ B_{\rho}(G)=B(G)$. Now apply \cite[Lemma 5.1]{Fendler-Lau}.

(b)$\Rightarrow$(c)$\Rightarrow$(d)$\Rightarrow$(e) follows from \cite[Theorem
5.3 and Theorem 6.1]{Fendler-Lau}.

(e)$\Longleftrightarrow$(a) follows from Theorem \ref{main.thm} and Proposition
\ref{converse}.

(b)$\Rightarrow$(f)$\Rightarrow$(g)$\Rightarrow$(h) follows from \cite[Theorem
5.3 and Theorem 6.1]{Fendler-Lau}.

(h)$\Rightarrow$(a) follows from \cite[Theorem 4.2 (iii)]{Miao}.

(i)$\Longleftrightarrow$(j)$\Longleftrightarrow$(k)$\Longleftrightarrow$(a) follow from \cite[Theorem 4.6]{Lau-Ulger}.

(a)$\Longrightarrow$(l) by amenability of $G$ and Peter-Weyl theorem \cite[Chapter 7]{Folland}.

(l)$\Longrightarrow$(a) follows from Theorem \ref{Bagget}.

$\blacksquare$

 A locally
compact Hausdorff topological space $X$ is scattered if $X$ does not
contain any non-empty perfect subset. Equivalently, every non-empty subset of
$X$ has an isolated point. By a theorem of Jensen \cite{Jensen 78}, $X$ is scattered
if and only if the C*-algebra $C_{0}(X)$ is scattered.
Some authors use the term {\it dispersed} instead of scattered. The following
example introduces some scattered spaces.

\begin{example} \label{Sctr}
  For a separable Hilbert space $\mathcal{H}$, the C*-algebra $\mathcal{K}(\mathcal{H})$
is scattered. Consider $\mathbb{N}$ with its usual topology and its one point compactification
  ${\mathbb{N}}^{*}=\mathbb{N} \cup \{\infty\}$. Both $\mathbb{N}$ and
  ${\mathbb{N}}^{*}$ are scattered. The set $\Omega=\{\frac{1}{n}: n\in \mathbb{N} \}\cup \{0\} $ of reciprocals of naturals with zero added
  provides another example of a compact scattered space. The set $\Omega$ has the nice property of being scattered but not discrete. We use it in our counter example.
\end{example}

Randrianantoanina \cite[Corrolary 3.8]{Randrian} has shown that if a scattered
C*-algebra $\mathcal{A}$ is a dual space, then the Banach space ${ \mathcal{A} }^{*}$ has weak* fpp for left reversible semigroups. Then, he asked \cite[Question
3.9]{Randrian} if this result is true without the dual assumption on ${
\mathcal{A} }$. In the following, we give a negative answer to this. As a result, the condition of being dual on $\mathcal{A}$ in \cite[Corrolary 3.8]{Randrian} cannot be removed even for abelian C*-algebras. We remark that an abelian C*-algebra $\mathcal{A}$ is a dual space if and only if its spectrum $\hat{A}$ is a hyperstonean space ( see \cite[p.109]{Takesaki} for more details).

\begin{thm} \label{Randr}
For $\Omega=\{\frac{1}{n}: n\in \mathbb{N} \}\cup \{0\}$, the C*-algebra $\mathcal{A}=C(\Omega)$ is scattered, but its dual Banach space ${ \mathcal{A} }^{*}$ does not have weak* fpp for left reversible semigroups.
\end{thm}
{\bf Proof.} 
Obviously by Jensen's theorem \cite{Jensen 78}, $\mathcal{A}$ is scattered since $\Omega=\hat{\mathcal{A}}$ is scattered. We claim that $\mathcal{A}$ is separable. Note that, $\Omega$ is a metric space under the absolute value of real numbers and $\{\frac{1}{n}: n\in \mathbb{N} \}$ is a dense subset of it. For each  $n\in \mathbb{N}$, put $f_{n}(x)=|x-\frac{1}{n}|$. By continuity of $|.|$, it is obvious that each $f_{n}$
is continuous. Also, the family $\{1\}\cup\{f_{n}: n\in \mathbb{N}\}$ of real-valued functions is a separating family in $\mathcal{A}$. By Stone-Weierstrass theorem, the algebra generated by $\{1\}\cup\{f_{n}: n\in \mathbb{N}\}$ is dense in $\mathcal{A}$. Using the rational coefficient for the linear combinations of  $\{1\}\cup\{f_{n}: n\in \mathbb{N}\}$, we see that $\mathcal{A}$ is separable. Suppose to the contrary that ${ \mathcal{A} }^{*}$ has weak* fpp for left reversible semigroups. By \cite[Theorem 2.1]{Fendler-Leinert}, $\hat{\mathcal{A}}=\Omega$ must be discrete. But, this is obviously false since for each $m\neq n$, $|\frac{1}{n}-\frac{1}{m}|\neq 0$. This contradiction shows that the answer to question 3.9 in \cite{Randrian} is negative. $\blacksquare$

\section{Bruck's Generalized Conjecture for preduals of von Neumann algebras}

Bruck \cite{Bruck} has shown that if a Banach space $E$
has weak fpp, then it has weak fpp for abelian semigroups.
We call the following statement {\it Bruck's Generalized Conjecture (BGC)}.

(BGC) If a Banach space $E$ has weak fpp, then it has weak fpp for any left
reversible semi-topological semigroup $S$.

It is a brilliant idea to pass from the private fixed points of single mappings
to the common fixed point of the whole mappings.
There are several techniques for this passage (for example,
asymptotic center property and lim-sup
property as discussed in \cite{Lim},\cite{Fendler-Lau} and\cite{Randrian}).
But, these techniques are interwoven with the
geometric and algebraic properties of the Banach space in question.
On the other hand, Lim's fixed point theorem
\cite[Theorem 3]{Lim-wk-fpp} imposes a strong condition ( called normal structure) on the
space in question to assure the existence of a common fixed point for left
reversible semi-topological semigroups. But, many
good spaces like $c_0$ lack normal structure and still have weak fpp. Meanwhile, amenability of groups has an equivalent
formulation in terms of fixed points (see\cite[p.61]{Pat}).  So, we
suspect that the BGC would be true in general unless one is dealing with well-behaved semigroups
or special spaces as in \cite{Bruck}, \cite{Fendler-Lau}, \cite{Lim}, and \cite{Randrian}. 

We are going to show that the BGC is true for the preduals of von Neumann algebras. Before doing so, we recall some terminology. Let $\mathcal{M}$ be a von Neumann algebra, ${\mathcal{M}}^{'}$ be its commutant and ${ \mathcal{M}}_{prj}$ denotes the set of all projections in $\mathcal{M}$. A projection $z\in { \mathcal{M}}_{prj}$ is called {\it central} if $z\in \mathcal{M}\cap {\mathcal{M}}^{'}$. A projection $e\in { \mathcal{M}}_{prj}$ is called {\it abelian} if $e{\mathcal{M}}e=\mathcal{M}$. The von Neumann algebra $\mathcal{M}$ is called to be of {\it Type I} if every non-zero central projection in it dominates a non-zero
abelian projection. $\mathcal{M}$ is called a {\it factor} if its center consists of all complex multiple of the identity, i.e., $\mathcal{M}\cap {\mathcal{M}}^{'}=\mathbb{C}\mathcal{I}$. Unfortunately, the Type I terminology
conflicts with the same terminology for C*-algebras since Type I von Neumann algebras need not be Type I C*-algebras in general. A {\it factor of Type I } is just a von Neumann algebra which is both a factor and a Type I
von Neumann algebra. Consult \cite{Takesaki} for more details on von Neumann algebras.

The proof of the following theorem is inspired by Randrianantoanina's ideas \cite{Randrian} and generalizes
Bruck's fixed point theorem \cite{Bruck} for the preduals of von Neumann algebras.

\begin{thm} \label{vn}
Let $\mathcal{M}$ be a von Neumann algebra with the predual ${\mathcal{M}}_{*}$. If ${\mathcal{M}}_{*}$ has weak fpp, then it has weak fpp for left reversible semigroups.
\end{thm}

{\bf Proof.} By \cite[Lemma 2.1]{Fendler-Leinert}, $\mathcal{M}$ is atomic, of Type I and has RNP.
 An application of \cite[p. 299]{Takesaki}, reveals that 
\begin{center}
$\mathcal{M}=l^{\infty}-\oplus_{i\in I} L^{\infty}(\mu_{i}) {\overline{\otimes}}\mathcal{B}({\mathcal{H}}_{i})$
\end{center}
and since $\mathcal{M}$ is atomic, every measure $\mu_{i}$ occurring in the above sum is atomic. So, without loss of generality, we write 
\begin{center}
$\mathcal{M}=l^{\infty}-\oplus_{i\in I} l^{\infty}(\mu_{i}) {\overline{\otimes}}\mathcal{B}({\mathcal{H}}_{i}).$
\end{center}
By uniqueness of preduals of von Neumann algebras, we have
 
\begin{center}
$\mathcal{M}_{*}=l^{1}-\oplus_{i\in I} l^{1}(\mu_{i}) {\overline{\otimes}}\mathcal{T}({\mathcal{H}}_{i})=
l^{1}-\oplus_{i\in I} l^{1}(\mu_{i} ,\mathcal{T}({\mathcal{H}}_{i})).$
\end{center}
Hence, by a diagonalization argument, $\mathcal{M}_{*}$ can be embedded isometrically into a closed subspace of $\mathcal{T}(l^{2}-\oplus_{i\in I} l^{2}(\mu_{i})\otimes {\mathcal{H}}_{i})$. By \cite[Theorem 3.6]{Randrian}, $\mathcal{T}(\mathcal{H})$, for any Hilbert space $\mathcal{H}$,  has weak fpp for left reversible semigroups. Hence, the result follows. $\blacksquare$

The following example shows one must deal with non-expansive actions in BGC. Recall that the non-amenable
group $SL_{2}(\mathbb{R})$ can act on the unit sphere $\mathbb{S}$ of
${\mathbb{R}}^3$ by the M\"{o}bius transformations. Besides, the unit sphere in ${\mathbb{R}}^3$ is the one point
compactification of the complex plane with south pole corresponds to zero and north pole corresponds to
$\infty$. 

\begin{example}
The action of the group $SL_{2}(\mathbb{R})=\{ \begin{pmatrix}
a&b\\
c&d
\end{pmatrix} : a, b, c, d \in \mathbb{R}; ad-bc=1 \}$ on $\mathbb{S}$ by M\"{o}bius
transformations is
separately continuous and fixed point free while each member of $SL_{2}(\mathbb{R})$ has
a fixed point.
\end{example}
{\bf Proof.} The action is given by
\begin{center}
\begin{align*}
 &SL_{2}(\mathbb{R})\times \mathbb{S}\longrightarrow \mathbb{S}\\
&{\begin{pmatrix}
a&b\\
c&d
\end{pmatrix}} \times {z}\longrightarrow \frac{az+b}{cz+d}
\end{align*}
\end{center}
Obviously, each member of the group defines a continuous ( not non-expansive, not affine)
mapping (called M\"{o}bius transformation) on the compact set $\mathbb{S}$.
So, the action is separately continuous.
By Schauder's fixed point theorem, the action of each member of the group has a fixed point.
 Also by the very formula of each transformation, one
can find the fixed points directly. By an easy calculation,
one can find two transformations without a common fixed point. Hence, the action does
not have a common fixed point.$\blacksquare$

The next example, uses a semigroup which is right reversible but not left reversible. As a result, BGC is not true for right reversible semigroups.

\begin{example}
Let $S=\{f, g\}$ be the discrete semigroup with the multiplication defined by $ab=a$ for all $a,b\in S$. Obviously, $S$ is not left reversible. Let $\alpha$ and $\beta$ be two real numbers with $0<\alpha<\beta<1$. Put $K=\{(a_n): (a_n)\in c_0(\mathbb{Z}) , \alpha\leq a_0\leq\beta; \quad a_{n}=0 \quad \text{for all} \quad n\neq 0\}$. By an easy observation, we see that $K$ is convex and weakly closed set in $c_0(\mathbb{Z})$. By James' theorem \cite{James}, one can check that $K$ is also weakly compact.  Now, consider the action $S\times K\longrightarrow K$
defined by $f((a_n))=(...,0,\alpha,0,...)$ and $g((a_n))=(...,0,\beta,0,...)$ where $\alpha$ and $\beta$ are placed at the zero location. Obviously, each member of $S$ defines a non-expansive mapping on $K$, so the action is separately continuous and non-expansive. On the other hand, $c_0(\mathbb{Z})$ has weak fpp by \cite[Corolarry 4.3]{Lau-Mah-Ulger}. Either way, $c_0(\mathbb{Z})=C^{*}(\mathbb{T})$ is a closed C* sub-algebra of $\mathcal{K}(L^{2}(\mathbb{T}))$ and the ideal of compact operators has weak fpp by \cite{Dowling} and \cite{Garcia}. So, $c_0(\mathbb{Z})$ has weak fpp. Hence, each function in $S$ has a fixed point on $K$. But, there is no common fixed point for the action. Remember that, $S$ is not left reversible.
\end{example}

We would like to conclude our work with the following questions:

{\bf Question 1.} Is there any counterexample of a non-expansive action to BGC?

{\bf Question 2.} Is BGC true for amenable semigroups?

The norm of C*-algebras are completely determined by the algebraic properties of the space. Inspired by Theorem \ref{vn}, we ask:

{\bf Question 3.} Is BGC true for non-unital C*-algebras?

From Theorem \ref{vn}, we know the preduals of von Neumann algebras enjoy BGC property. Now, we ask:

{\bf Question 4.} Can we characterize all Banach spaces with BGC property?

{\bf Question 5.} Can one prove Theorem \ref{main.thm} for weak* fpp without using the second
countability assumption on $G$?

{\bf Acknowledgment.} The author is much indebted to Professor Gero Fendler and would like to thank him for his invaluable suggestions and remarks. He would also like to thank Professor Anthony To-Ming Lau for his instructions and hints which made this paper better.
 \\


\end{document}